\documentclass[11pt,a4paper]{amsart}
\usepackage{amssymb,amsmath,epsfig,graphics,mathrsfs,enumerate,verbatim}\usepackage{bm}
\usepackage[pagebackref,colorlinks=true,linkcolor=blue,citecolor=blue]{hyperref}
\usepackage{fancyhdr}
\usepackage{hyperref}
\hypersetup{
 colorlinks   = true,
 urlcolor     = blue,
 linkcolor    = blue,
 citecolor   = red ,
 bookmarksopen=true
}

\usepackage{amsmath}
\usepackage{amsfonts}
\usepackage{amssymb}
\usepackage{amsthm}
\usepackage{epsfig,graphics,mathrsfs}
\usepackage{graphicx}
\usepackage{dsfont}

\usepackage[usenames, dvipsnames]{color}

\usepackage{hyperref}

\def \N {\mathbb{N}}

\def \phi {\varphi}

\def \R {\mathbb{R}}

\newcommand{ \G}{\Gamma}

\newcommand{\Ba}{\mathscr B_a}

\newcommand{\Rn}{\mathbb R^n}

\newcommand{\Om}{\Omega}
\newcommand{\bomega}{{\partial\Omega}}

\newcommand{\p}{\partial}
\newcommand{\IO}{\int_\Omega}
\newcommand{\bG}{\mathbb {G}}

\newcommand{\la}{\lambda}

\numberwithin{equation}{section}

\newcommand{\beq}{\begin{equation}}
\newcommand{\bea}[1]{\begin{array}{#1} }
\newcommand{\eeq}{ \end{equation}}
\newcommand{\ea}{ \end{array}}

\newcommand{\ve}{\varepsilon}

\newcommand{\tr}{\operatorname{tr}}
\newcommand{\divg}{\operatorname{div}}

\newcommand{\nh}{\nabla_H}

\newcommand{\sa}{\langle}
\newcommand{\da}{\rangle}
\newcommand{\gr}{\nabla\,}
\newcommand{\vola}{\, r^adx}
\newcommand{\surfa}{\, r^a d\sigma}

\newcommand{\Bab}{\mathscr B_{\boldsymbol{a}}}
\newcommand{\Bai}{\mathscr B_{a_i}}
\newcommand{\Baj}{\mathscr B_{a_j}}
\newtheorem{theorem}{Theorem}[section]
\newtheorem{thrm}{Theorem}[section]
\newtheorem{lemma}[theorem]{Lemma}
\newtheorem{proposition}[theorem]{Proposition}
\newtheorem{prop}[theorem]{Proposition}

\newtheorem{cor}[theorem]{Corollary}
\newtheorem{remark}[theorem]{Remark}

\numberwithin{equation}{section}

\begin{document}

\title[Overdetermined, etc.]{Overdetermined fractional Serrin problem}



\date{}

\begin{abstract}
We solve a version of the Serrin overdetermined problem for the Weinstein operator involving a Bessel operator.
\end{abstract}

\author{Nicola Garofalo}
\address{School of Mathematical and Statistical Sciences\\ Arizona State University}\email[Nicola Garofalo]{nicola.garofalo@asu.edu}

\thanks{N. Garofalo's present affiliation is the Arizona State University, but his work on this project started when he was on the faculty at the University of Padova. He was supported in part by a Progetto SID (Investimento Strategico di Dipartimento): ``Aspects of nonlocal operators via fine properties of heat kernels", University of Padova (2022); and by a PRIN (Progetto di Ricerca di Rilevante Interesse Nazionale) (2022): ``Variational and analytical aspects of geometric PDEs". He was also partially supported by a Visiting Professorship at the Arizona State University}

\author{Dimiter Vassilev}

\address{University of New Mexico\\
Department of Mathematics and Statistics\\
311 Terrace Street NE\\
Albuquerque, NM  87106}
\vskip 0.2in
\email{vassilev@unm.edu}

\maketitle

\tableofcontents

\section{Introduction}\label{S:intro}

Overdetermined boundary value problems have long served as a powerful source of rigidity phenomena in both geometric and analytic contexts. A classical example in this vein is the celebrated theorem of Serrin \cite{Se}, which asserts that if a function $u:\Om\to \R$
solves the torsion problem
\[
\Delta u = - 1\ \ \text{in}\ \ \Om,\ \ \ \ u = 0\ \ \text{on}\ \ \p \Om,
\]
for a bounded, connected domain $\Om\subset \Rn$, with $C^2$ boundary, then the overdetermined condition $|\nabla u| = c$ on $\p \Om$ holds if and only if $\Om$ is a ball $B(x_0,R)$. In that case $u$ is explicitly given by
\[
u(x) = \frac{R^2 - |x-x_0|^2}{2n}.
\]
As it is well-known, in his paper \cite{We} Weinberger gave an alternative proof of Serrin's result based on an ingenious combination of the strong maximum principle and integral identities. His approach was subsequently generalized by Lewis and the first named author in \cite{GL} to quasilinear operators modelled on the $p$-Laplacian.

In this work, we establish a Serrin-type symmetry result for a class of degenerate elliptic operators that arise naturally in the theory of generalized axially symmetric potentials. As we subsequently explain, our results find applications in some overdetermined problems in sub-Riemannian geometry.
For $a\ge 0$ we consider the \emph{Weinstein operator}
\begin{equation}\label{L}
L_a u = \p_{rr} u + \frac ar \p_r u + \Delta_{y} u,\ \ \ \ \ (r,y)\in \R^+\times \R^k,
\end{equation}
which couples a Bessel-type operator in the variable $r$,
\[
\Ba = \p_{rr}  + \frac ar \p_r = r^{-a} \p_r(r^a \p_r),
\]
with the standard Laplacian in the tangential variables $y\in \R^k$. When $a\in \N$ the operator $L_a$ corresponds to the Laplacian in $\R^{a+1+k}$ acting on functions which are spherically symmetric in the first $a+1$ variables. However, for arbitrary positive $a>0$ this geometric interpretation breaks down, and the analysis becomes significantly more delicate due to the degeneracy of the operator at $r=0$.

The second order operator \eqref{L} was considered by  Weinstein in \cite{We48, We2} in his study of generalized
axially symmetric potentials. With different objectives, it was subsequently studied in two important papers by Muckenhoupt and Stein \cite{MS} and Talenti \cite{Ta}. Moreover,  $L_a$ plays a key role in the celebrated extension procedure of Caffarelli and Silvestre for the fractional powers of the Laplacian $(-\Delta)^s$, see \cite{CS}. In particular, for $s\in (0,1)$, if one sets $a = 1-2s\in (-1,1)$, and solves the Dirichlet problem in $\R^+\times \R^k$
\[
L_a U = 0,\ \ \ \ \ \ \ \ \ U(0,y) = u(y),
\]
then the extension function $U(r,y)$ satisfies the  Dirichlet-to-Neumann condition
\[
- \frac{2^{2s-1} \G(s)}{\G(1-s)} \underset{r\to 0^+}{\lim} r^{1-2s} \p_r U(r,y) = (-\Delta)^s u(y),
\]
which connects the nonlocal operator $(-\Delta)^s$ in $\R^k$ to the local degenerate elliptic equation $L_a U = 0$ in one higher dimension.

Henceforth, we shall routinely denote $\R^{k+1} = \R_r\times \R^k_y$, and $\R^{k+1}_+ = \R^+_r\times \R^k_y$. We also indicate by $\Omega^\star\subset \R^{k+1}$ a bounded, connected open set, symmetric with respect to the hyperplane  $r=0$, and denote by $\Om=\Om^\star \cap \R^{k+1}_+$. Our main result is the following rigidity theorem for the overdetermined problem associated with the operator \eqref{L}.

\begin{thrm}\label{t:serrin}
Let  $\Omega^\star\subset \R^{k+1}$ be as stated, and assume additionally that it be piecewise $C^1$. The overdetermined problem
 \begin{equation}\label{e:serrin}
L_a u = - 1\quad \text{in}\ \Om^\star\setminus\{r=0\},
\qquad\quad
u = 0\quad  \text{on}\ \p \Om^\star,\qquad |\nabla u| = c  \quad \text{on}\  \p\Om^\star,
\end{equation}
admits a solution  $u\in C^2(\Om^\star\setminus \{r=0\})\cap C^1(\bar\Om^\star)$  if and only if for some $R>0$ and $y_0\in \R^k$, one has
\[
\Om^\star = B_R((0,y_0))=\{(r,y)\in \R^{k+1}\mid r^2 +|y - y_0|^2<R^2\},
\]
and $$u(r,y)=\frac {R^2-r^2- |y - y_0|^2}{2(a+1+k)}.$$
\end{thrm}
The proof of Theorem \ref{t:serrin} is inspired to the above mentioned Weinberger's alternative approach to Serrin's result. Our investigation was primarily motivated by our recent work on overdetermined problems in sub-Riemannian geometry. In the linear case $p=2$ of \cite[Theor. 1.7]{GVjfa}, we proved that in a group of Heisenberg type $\bG$ (with horizontal gradient $\nh$ and Laplacian $\Delta_H$), satisfying a certain structural transitivity invariance, which we called \emph{Property $(H)$}, if a domain with partial symmetry $\Om\subset \bG$ supports a solution $f$ to the overdetermined problem
\begin{equation}\label{i3}
\begin{cases}
\Delta_{H} f = - |z|^2\ \ \ \ \ \ \ \ \ \ \ \ \ \ \ \ \ \ \ \ \ \ \ \ \ \ \ \ \ \text{in}\ \Om,
\\
f = 0,\ \ \ \ |\nh f|_{\big|\p\Om} = c |z|, \ \ \ \ \ \ \ \ \ \ \ \ \ \text{on}\ \p \Om,
\end{cases}
\end{equation}
then there exist $R>0$ and $\sigma_0\in \R^k$ such that $\Om$ is a ball of the Koranyi gauge centred at $(0,\sigma_0)$ with radius $R$, i.e.
\[
\Om = B_R(0,\sigma_0) = \{(z,\sigma)\in \bG\mid |z|^4 + 16 |\sigma-\sigma_0|^2 < R^4\},
\]
and $f$ is explicitly determined.
Using the Property (H) and a suitable change of variables, we showed that, remarkably, the sub-Riemannian problem \eqref{i3} can be reduced to the Weinstein-type equation considered in the present work, with the Bessel parameter given by $a = \frac m2-1$, where $m$ represents the dimension of the horizontal layer of the Lie algebra of $\bG$. This connection underscores the broader significance of Theorem \ref{t:serrin}.

\subsection{Outline of the paper}
A brief description of the organization of the paper is as follows. In Section \ref{S:CDI} we develop a version of the Bakry-Emery Gamma calculus for the singular operator \eqref{L}, which we subsequently apply to  a suitable $P$-function \`a la Weinberger $P_a$, which we introduce in \eqref{Pa}. We show in Corollary \ref{T:CDLa} that $L_a$ verifies the Bakry-Emery \emph{curvature-dimension inequality} $\operatorname{CD}(0,a+1+k)$, with equality in \eqref{P} characterizing quadratic solutions of the form
\begin{equation}\label{e:cdi eq for winstein}
u(r,{y}) = \alpha(r^2+|{y} - {y}_0|^2) +\gamma,
\end{equation}
for some ${y}_0\in\R^k$ and $\alpha, \gamma \in \R$. This is then used to prove the key result of the section, Theorem \ref{T:P}.

In Section \ref{s:max principle} we generalize to \eqref{L} the two-dimensional strong maximum principle of Muckenhoupt and Stein \cite[Theor. 1]{MS}.

Section \ref{S:main} contains several technical results critical to the proof of Theorem \ref{t:serrin}. These include:
\begin{itemize}
\item A regularity result across the singular set ${r = 0}$ (Proposition \ref{p: normal der vanish}).
\item A Rellich-type identity (Proposition \ref{p:Pohozaev}),
\item An integral identity for the function $P_a$ (Theorem \ref{t:P integral for serrin}).
\end{itemize}
Combining these ingredients, we complete the proof of the main Theorem \ref{t:serrin} in Section \ref{S:main}, by showing that the
function $P_a$ must be constant throughout $\Om^\star$.

A final comment is in order. In Theorem \ref{t:serrin} we do not assume a priori that the solution $u$ is symmetric with respect to $r$. This property is in fact deduced as a consequence of the above discussed maximum principle for the operator $L_a$. On the other hand, we do assume that $u\in C^1(\bar \Om^\star)$. Had we assumed a priori evenness in $r$, we could have invoked the interesting \cite[Theorem 1.1]{STV}, which would have guaranteed $u\in C^\infty(\Om^\star)$, and therefore in particular across the singular hyperplane $\{r = 0\}$.

%%%%%%%%%%%%%%%%%%%%%%%%%%%%%%%%%%%%%%%%%%%%%%%%%%%%%%%%

\section{Curvature-dimension inequalities}\label{S:CDI}

In this section we develop a version of the Bakry-Emery Gamma calculus (see \cite{BE} and \cite{BGL}) for the singular operator \eqref{L}, and prove for the latter an auxiliary result, Theorem \ref{L:nicetoo}, which we then use

We start with a simple consequence of the Cauchy-Schwarz inequality which will be useful in the sequel.

\begin{lemma}\label{l:aCS}
Given numbers $w_i>0$ and $A_i\in\R$,  $i=1,\dots,n$, one has the inequality
\begin{equation}\label{e:elem inequality}
\sum_{i=1}^n\frac {1}{w_i}A_i^2\geq \frac{1}{\sum_{i=1}^n w_i}\left( \sum_{i=1}^n A_i\right)^2.
\end{equation}
Furthermore, equality holds if and only  $A_i=\lambda w_i$ for $i=1,\dots,n$, for some $\lambda\in\R$.
\end{lemma}

We next recall some basic definitions. Given a diffusion operator $L$ symmetric with respect to a smooth measure $\mu$ in $\Rn$, the carr\'e du champ associated with $L$ is defined by the equation
\begin{equation}\label{GammaD}
\G^L(u,v) = \frac 12 \left[L(u v) - u Lv - v Lu\right].
\end{equation}
It is obvious that $\G(u,v) = \G(v,u)$, and one easily verifies that
\begin{equation}\label{ui}
\G(u_1+...+u_k,v) = \sum_{i=1}^k \G(u_i,v).
\end{equation}
The iterated carr\'e du champ is given by
\begin{equation}\label{Gamma2D}
{\G_2^L}(u,v) = \frac 12 \left[L \G^L(u, v) - \G^L(u,Lv) - \G^L(v,Lu)\right].
\end{equation}
If $u = v$, we simply write $\G^L(u) = \G^L(u,u)$ and $\G^L_2(u) = {\G_2^L}(u,u)$.

When $L = \Delta$ is the Laplacian in $\R^k$, one easily verifies from \eqref{GammaD} and \eqref{Gamma2D} that
\begin{equation}\label{1}
\G^{\Delta}(u) = |\nabla u|^2,\ \ \ \ \ \ {\G_2^{\Delta}}(u) = ||\nabla^2 u||^2,
\end{equation}
where we have indicate with $\nabla^2$ the Hessian matrix of $u$ with respect to the canonical basis of $\R^k$.
The important inequality
\begin{equation}\label{cs}
{\G_2^{\Delta}}(u) \ge \frac{1}{k} (\Delta u)^2,
\end{equation}
referred to as the  \emph{curvature-dimension} inequality CD$(0,k)$,
is a well-known consequence of \eqref{Gamma2D}, \eqref{1}, the Bochner identity
\[
\Delta(|\nabla u|^2) = 2 ||\nabla^2 u||^2 + 2\sa \nabla u,\nabla(\Delta u)\da,
\]
and of the Newton inequality for symmetric matrices $A\in \mathbb M_{k\times k}(\R)$
\[
||A||^2 \ge \frac 1k (\tr(A))^2.
\]

More in general, we might consider the following situation. Let $M_1,...,M_n$ be Riemannian manifolds, with assigned diffusion operators $L_1,...,L_n$, and invariant measures $\mu_1,...,\mu_n$. On the product manifold $M = M_1\times ... \times M_n$, we consider the operator $L = L_1+...+L_n$, with $\mu = \mu_1\otimes ... \otimes \mu_n$. Then it is easy to verify that for every two functions $u, v\in C^\infty(M)$, one has
\begin{equation}\label{sum}
\G^L(u,v) = \sum_{i=1}^n \G^{L_i}(u,v),
\end{equation}
where $\G^{L_i}(u,v)$ has the obvious meaning that the differential operator $L_i$ acts only on the variable $x_i\in M_i$. Next, we claim that for the iterated \emph{carr\'e du champ} one has
\begin{equation}\label{D}
\G^L_2(u,v) = \sum_{i=1}^n \G_2^{L_i}(u,v)+ \frac 12 J(u,v),
\end{equation}
where we have set
\begin{equation}\label{JL}
J_L(u,v) = \sum_{i=1}^n \sum_{j\not= i} \left\{L_i\G^{L_j}(u,v) -  \G^{L_i}(u,L_j v) - \G^{L_i}(v,L_j u)\right\}.
\end{equation}
This follows by the ensuing computation. For $u\in C^\infty(M)$ one has from \eqref{ui}, \eqref{sum}
\begin{align*}
&\G^L_2(u,v) = \frac 12\left[L \G^L(u,v)  - \G^L(u,Lv)- \G^L(v,Lu)\right]
\\
& = \frac 12\bigg[\sum_i L_i \G^{L_i}(u,v) + \sum_i \sum_{j\not= i} L_i\G^{L_j}(u,v) - \sum_i \G^{L_i}(u,L_i v) - \sum_i \sum_{j\not= i} \G^{L_i}(u,L_j v)
\\
& - \sum_i \G^{L_i}(v,L_i u) - \sum_i \sum_{j\not= i} \G^{L_i}(v,L_j u)\bigg]
\\
& = \sum_{i} \G_2^{L_i}(u,v) + \frac 12 J_L(u,v),
\end{align*}
which proves \eqref{D}. In particular, we find from \eqref{JL}
\begin{equation}\label{JLu}
J_L(u,u) = \sum_{i=1}^n \sum_{j\not= i} \left\{L_i\G^{L_j}(u) - 2 \G^{L_i}(u,L_j u)\right\}.
\end{equation}
These considerations lead to the following result.
\begin{proposition}\label{P:prod}
Suppose that the triple $(M_i,\mu_i,L_i)$ satisfy that curvature-dimension inequality CD$(\rho_i,k_i)$ for $i=1,...,n$, i.e.,
\begin{equation}\label{Li}
\G^{L_i}_2(u) \ge \rho_i \G^{L_i}(u) + \frac{1}{k_i} (L_i u)^2
\end{equation}
for every $u\in C^\infty(M_i)$. Then the triple $(M,\mu,D)$ satisfies the condition
\begin{equation}\label{functD}
\G_2^L(u) \ge \rho \G^L(u) + \frac 1k (Lu)^2+ \frac 12 J_L(u,u),
\end{equation}
where $\rho = \min\{\rho_1, ..., \rho_n\}$, $k = k_1+...+k_n$.
\end{proposition}

\begin{proof}
To prove \eqref{functD}, we use \eqref{D} and \eqref{Li}, which give
\begin{align*}
\G^L_2(u,v) & \ge \sum_{i=1}^n\left(\rho_i \G^{L_i}(u) + \frac{1}{k_i} (L_i u)^2\right)+ \frac 12 J_L(u,u)
\\
& \ge \rho \G^L(u) + \sum_{i=1}^n \frac{1}{k_i} (L_i u)^2 + \frac 12 J_L(u,u)
\\
& \ge \rho \G^L(u) + \frac{1}{\sum_{i=1}^n k_i} (L u)^2+ \frac 12 J_L(u,u),
\end{align*}
where in the last inequality we have applied \eqref{e:elem inequality}
in Lemma \ref{l:aCS}.

\end{proof}

\subsection{Curvature-dimension inequality for  Bessel operators}\label{ss: CDI Bessel}
Given a number $a\ge 0$, consider the Bessel operator \eqref{Ba} in $\R^+$
\begin{equation}\label{Ba}
\Ba = \p_{rr}  + \frac ar \p_r = r^{-a} \p_r(r^a \p_r),
\end{equation}
with its invariant measure $d\mu = r^a dr$.
It was observed in \cite{G} that
\begin{equation}\label{2}
\G^{\Ba}(u) = (\p_r u)^2,\qquad  \G_2^{\Ba}(u) = (\p_{rr}u)^2 + \frac{a}{r^2} (\p_r u)^2.
\end{equation}
Furthermore, the following basic fact holds.

\begin{lemma}\label{L:nice}
For every $u\in C^\infty(\R^+)$ one has
\begin{equation}\label{nice}
{\G_2^{\Ba}} (u) \ge \frac{1}{a+1} (\Ba u)^2.
\end{equation}
If $a>0$, equality holds in \eqref{nice} if and only if $ u=\alpha r^2+\gamma$, with $\alpha,  \gamma\in \R$. When $a=0$ equality in \eqref{nice} holds trivially on any $u\in C^\infty(\R^+)$.
\end{lemma}

\begin{proof}
The inequality \eqref{nice} was noted in \cite[p. 111] {BakryQian00} and \cite[Prop. 3.1]{G} where it was also observed that \eqref{nice} is false (in fact, it gets reversed) when $a<0$. For completeness we present the elementary proof here. Applying Lemma \ref{l:aCS} with $w_1 = 1, w_2 = a$, and $A_1 = \p_{rr} u, A_2 = \frac ar \p_r u$, and noting that $\Ba u = A_1 + A_2$, we obtain
\begin{align*}
& \frac{1}{a+1} (\Ba u)^2 \le (\p_{rr} u)^2 + \frac 1a (\frac ar \p_r u)^2 = (\p_{rr} u)^2 + \frac a{r^2} (\p_r u)^2 = \G_2^{\Ba}(u),
\end{align*}
where in the last equality we have used the second equation in \eqref{2}. This immediately shows that, when $a = 0$, equality holds for any $u\in C^\infty(\R^+)$. If instead $a>0$, from the case of equality in Lemma \ref{l:aCS}, we obtain at any $r>0$
\[
\p_{rr} u = \la(r),\ \ \ \ \ \frac ar \p_r u = \la(r) a,
\]
for some $\la(r)\in \R$. This gives at every $r>0$
\[
\p_r(r^{-1}\p_r u) = \frac 1r(\p_{rr} u - \frac 1r \p_r u) = 0,
\]
therefore there exists $\alpha\in \R$ such that $\p_r u = 2 \alpha r$, and then $u(r) = \alpha r^2 + \gamma$.

\end{proof}

In the language of Bakry-Emery, the inequality \eqref{nice} shows that the Bessel operator $\Ba$ satisfies the curvature-dimension inequality CD$(0,a+1)$ in $\R^+$. We next prove a generalization of Lemma \ref{L:nice} which is of independent interest. Given numbers $a_1,...,a_n\ge 0$, we let $\boldsymbol{a} = (a_1,...,a_n)$. Consider the sum of the $n$ Bessel operators $\Bai$,
\begin{equation}\label{e:sum of n bessel}
\Bab = \sum_{i=1}^n \Bai = \sum_{i=1}^n\left(\p_{r_i r_i} + \frac{a_i}{r_i} \p_{r_i}\right),
\end{equation}
acting on the space $\mathcal{R}=\prod_{i=1}^n \mathcal R_i$, where
$$\mathcal R_i= \begin{cases}\label{e:domain for sum of Bessel}
    \R^+  & \text{if $a_i>0$;} \\
    \R , & \text{if $a_i=0$.}
  \end{cases}
$$

\begin{thrm}\label{L:nicetoo}
Let $n\geq 2$. For every $u\in C^\infty(\mathcal{R})$ the following curvature-dimension inequality holds
\begin{equation}\label{e:CDI n Bessel}
\Gamma_2^{\Bab}(u)
\geq  \frac {1}{n+\sum_{i=1}^n a_i }  \Bab (u)^2.
\end{equation}
Equality holds in \eqref{e:CDI n Bessel} if and only if
\begin{equation}\label{e:CDI equality Bessels}
u=\gamma+ \alpha \sum_{i=1}^n \left(r_i^2+ \beta_i r_i\right),
\end{equation}
where $\beta_i=0$ if $a_i>0$.
\end{thrm}

\begin{proof}
If we apply \eqref{functD} to $M_i = \R^+$, with $L_i = \Bai$, $\rho_i = 0$ and $k_i = a_i+1$, then the inequality \eqref{e:CDI n Bessel} immediately follows from \eqref{functD} and from \eqref{nice} in Lemma \ref{L:nice}, once we prove that for $L = \Bab$ as in \eqref{e:sum of n bessel}, we have
\begin{equation}\label{B}
J_L(u,u) \ge 0.
\end{equation}
From \eqref{JL} we find
\begin{align*}
J_L(u,u) & = \sum_{i=1}^n \sum_{j\not= i} \left\{\Bai\G^{\Baj}(u) - 2 \G^{\Bai}(u,\Baj u)\right\}.
\end{align*}
Keeping \eqref{2} in mind, we have
\[
\G^{\Baj}(u) = (\p_{r_j} u)^2.
\]
Next, given a function $v$, we have
\[
\Bai(v^2) = 2 v \Bai v + 2 (\p_{r_i} v)^2.
\]
Combining these equations, we thus find
\begin{align*}
\Bai\G^{\Baj}(u) & = \Bai((\p_{r_j} u)^2) = 2 \p_{r_j} u \Bai(\p_{r_j} u) + 2 (\p_{r_i} \p_{r_j} u)^2.
\end{align*}
On the other hand, \eqref{2} again gives for $j\not= i$
\[
\G^{\Bai}(u,\Baj u) = \p_{r_i} u \p_{r_i}(\Baj u) = \p_{r_i} u \Baj(\p_{r_i} u).
\]
We finally obtain
\begin{align*}
J_L(u,u)  & = \sum_{i=1}^n \sum_{j\not= i}\left\{2 \p_{r_j} u \Bai(\p_{r_j} u) + 2 (\p_{r_i} \p_{r_j} u)^2 - 2
\p_{r_i} u \Baj(\p_{r_i} u)\right\}
\\
& = 2 \sum_{i=1}^n \sum_{j\not= i} (\p_{r_i} \p_{r_j} u)^2 \ge 0.
\end{align*}
This proves \eqref{B}, and therefore \eqref{e:CDI n Bessel}.
The case of equality in \eqref{e:CDI n Bessel} follows from the fact that we must have equality in all used inequalities, hence we have
\begin{align}
\p_{r_i} \p_{r_j} u & =0, \quad i=1,...n,, j\not=i,\label{e:3bn} \\
a_i\left(\p_{r_i r_i} u - \frac 1{r_i} \p_{r_i} u\right)  & =0, \label{e:1bn}\\
 \frac {\mathscr B_{a_1} u}{a_1+1} &=\dots =\frac {\mathscr B_{a_n}u}{a_n+1} \label{e:2bn} .
\end{align}
The vanishing of the mixed derivative \eqref{e:3bn} shows that
$\p_{r_j} u$ is independent of  the remaining variables, i.e, $\p_{r_j} u =f_j(r_j)$. We thus have
$$\frac {\mathscr B_{a_j} u}{a_j+1}  = \p_{r_j r_j}u =f'_j(r_j), \quad j=1,\dots,n.$$ It follows that for some constant $\lambda$ the Hessian $\nabla^2 u= \lambda I_n$, where $I_n$ is the identity operator in $\R^n$. Therefore, equality holds if and only if $u$ takes the claimed form \eqref{e:CDI equality Bessels}.

\end{proof}

Returning to the Weinstein operator $L_a$, we note that, with $n = k+1$, it can be itself considered as a sum of Bessel operators $\Bab$ as in \eqref{e:sum of n bessel}, where $\boldsymbol a = (a,0,...,0)$. If  we thus specialize Theorem \ref{L:nicetoo}  to $L_a$, we obtain the following result.

\begin{cor}\label{T:CDLa}
For any function $u\in C^\infty(\R^+\times \R^k)$ one has
\begin{equation}\label{P}
\G_2^{L_a}(u)\ \ge\ \frac{1}{a+1+k}\ (L_a u)^2.
\end{equation}
Equality occurs in \eqref{P} if and only if for some ${y}_0\in\R^k$ and $\alpha, \gamma \in \R$ we have
\begin{equation}\label{e:cdi eq for winstein}
u(r,{y}) = \alpha(r^2+|{y} - {y}_0|^2) +\gamma.
\end{equation}
\end{cor}

We explicitly note that \eqref{P} expresses the fact that the operator $L_a$ verifies the Bakry-Emery \emph{curvature-dimension inequality} $\operatorname{CD}(0,a+1+k)$. Corollary \ref{T:CDLa} will be used in the proof of Theorem \ref{T:P} below.

\subsection{The function $P_a$}\label{ss:P-function proof}

We continue the analysis by proving a key result concerning a (fractal) generalization of Weinberger's $P$-function to the setting of Theorem \ref{t:serrin}.

\begin{theorem}\label{T:P}
Consider an open set $\Om\subset \R^{k+1}_+$ and suppose that the function $u\in C^3(\Om)$ be a solution in $\Om$ of the equation
\begin{equation}\label{torsion}
L_a u = -1.
\end{equation}
Then, the function
\begin{equation}\label{Pa}
P_a \overset{def}{=} \G_{L_a}(u) + \frac 2{a+1+k}\ u
\end{equation}
satisfies the differential inequality in $\Om$
\begin{equation}\label{wow}
L_a P_a \ge 0.
\end{equation}
Finally, for a solution of \eqref{torsion}, equality is true in \eqref{wow} if and only if there exists ${y}_0\in\R^k$ and $\gamma \in \R$, such that for any $(r,y)\in \Om$,
\begin{equation}\label{e:serrin equality}
u(r,{y}) = \gamma - \frac{r^2+|{y}-{y}_0|^2}{2(a+1+k)}.
\end{equation}
\end{theorem}

\begin{proof}[Proof of Theorem \ref{T:P}]
From \eqref{torsion} and \eqref{Pa}, we have
\[
L_a P_a = L_a(\G_{L_a}(u)) - \frac 2{a+1+k}.
\]
Applying \eqref{Gamma2D} with  $L = L_a$, and noting that \eqref{torsion} gives
\[
\G_{L_a}(u,L_a u) = - \G_{L_a}(u,1) = 0,
\]
we find
\[
 {\G_2^{L_a}}(u) = \frac 12 L_a \left( \G_{L_a}(u)\right) - \G_{L_a}(u,L_a u) = \frac 12 L_a \G_{L_a}(u),
\]
Substituting the latter identity in the previous one, we obtain
\[
L_a P_a = 2 \left[ {\G_2^{L_a}}(u) - \frac 1{a+1+k}\right] \ \ge \ 0,
\]
where the (crucial) last inequality is justified by  Corollary  \ref{T:CDLa}.
For the second part of the theorem, if equality holds in \eqref{wow}, then we must have equality in \eqref{P}, and therefore $u$ takes the form \eqref{e:cdi eq for winstein} in Corollary   \ref{T:CDLa}. Finally,  \eqref{torsion} forces $\alpha = - \frac 2{a+1+k}$ in \eqref{e:cdi eq for winstein}, therefore \eqref{e:serrin equality} does hold.

\end{proof}

\section{The strong maximum principle}\label{s:max principle}

By the classical maximum principle of E. Hopf  \cite{Hopf}, given a connected open set $D\subset \R^{k+1}$ not intersecting the singular plane  $x_n = 0$, a function $u\in C^2(D)$ such that  $L_a u\geq 0$ cannot attain a maximum at an interior point, unless it is a constant. In this section we prove a strong maximum principle  for $L_au \geq 0$ in a domain $\Om^\star$ intersecting the singular plane, see Theorem \ref{T:smp}. In the two-dimensional case $(k=1)$, such strong maximum principle was first proved by Muckenhoupt and Stein in \cite[Theor. 1]{MS}. Their approach was based on anisotropic weighted spherical averaging, and we generalize such approach to the higher-dimensional case. This is of course very classical in spirit and in this respect we mention that mean-value formulas associated to the Weinstein operator were also found in \cite{We48, Kapilevich60, Weinacht62, Weinacht68, Kipriyanov97, Shishkina20}, see also see \cite{Delsarte38, Levitan49, Levitan51, Kingman, NahiaSalem94} for related results.
The generalized translations,
which commute with the Bessel operator,  are a starting point of the harmonic analysis related to the Weinstein operator and the associated hypergroup structure, see  \cite{Jewett75, BloomHeyer95, Trimeche01, ShishSit00}. Furthermore, the mean-value operators defined using the generalized translations intertwine with the operator $L_a$.

In what follows, we denote ${x}= (r,{y})\in \R^{k+1}$, and consider the Euclidean distance function in $\R^{k+1}$
\begin{equation}\label{d}
\rho (x)= \rho(r,{y}) = (r^2+|{y}|^2)^{1/2}=|x|.
\end{equation}
One easily checks that $\rho$ satisfies the following \emph{eikonal equation} in $\R^{k+1}\setminus\{(0,0)\}$,
\begin{equation}\label{eiko}
\G_{L_a}(\rho) = |\nabla \rho|^2 = (\p_r \rho)^2 + |\nabla_y \rho|^2 = 1.
\end{equation}
If $f\in C^2([0,\infty))$ and $u = f(\rho)$, then we easily find
\[
\p_r u = f'(\rho) \p_r \rho = \frac{f'(\rho)}{\rho} r , \ \ \ \ \ \nabla_y u = f'(\rho) \nabla_y \rho = \frac{f'(\rho)}{\rho} y.
\]
By elementary computations, this gives
\[
|\nabla u|^2  = f'(\rho)^2,\ \ \Ba u = f''(\rho) \frac{r^2}{\rho^2} + \frac{f'(\rho)}{\rho}(a+1 - \frac{r^2}{\rho^2}),\ \ \Delta_y u = f''(\rho) \frac{|y|^2}{\rho^2} + \frac{f'(\rho)}{\rho}(k - \frac{|y|^2}{\rho^2}).
\]
We thus find
\begin{equation}\label{radial0}
L_a u = \Ba u + \Delta_y u = f''(\rho) + \frac{a+k}{\rho} f'(\rho).
\end{equation}
With \eqref{radial0} in hand, we can prove that for every ${y}_0\in \R^k$ the function
\[
E(r,{y}) = - \frac{1}{(a+k-1){\sigma}_{a,k}} \rho(r,{y}-{y}_0)^{1-a-k},
\]
where ${\sigma}_{a,k}$ is given by \eqref{anivol} below, is a fundamental solution of $L_a$ with pole in $(0,{y}_0)$.

Consider now the Euclidean dilations $\delta_\la(x) = \la x$ in $\R^{k+1}$ and denote by $Z f = \frac{d}{d\la} (f\circ \delta_\la)\big|{_{\la =1}}$ their infinitesimal generator acting on a function $f$. One easily verifies that
\begin{equation}\label{e:radial field}
Z f=  \sa{x},\nabla f\da=r f_r+ \sa{y},\nabla_{y} f\da.
\end{equation}
Observe that for every $a\ge 0$ we have for $r\not= 0$
\begin{equation}\label{Za}
\operatorname{div}(|r|^a Z) = (a+1+k) |r|^a.
\end{equation}
Consider the Euclidean balls centred at $(0,{y}_0)\in \R^{k+1}$ with radius $t>0$
\[
B_t(0,{y}_0) = \{(r,{y})\in \R^{k+1}\mid \rho(r,{y}-{y}_0) < t\}.
\]
When ${y}_0 = 0$, we simply write $B_t$ for $B_t(0,0)$. We define the anisotropic volume and surface measure of $B_1$ as
\begin{equation}\label{anivol}
\omega_{a,k} = \int_{B_1} |r|^a dx,\ \ \ \ \ \ \sigma_{a,k} = \int_{\p B_1} |r|^a d\sigma,
\end{equation}
where we have denoted by $d\sigma$ the standard surface measure on $\p B_1$. We have the following.
\begin{lemma}\label{L:aniball}
For every $a>0$ and $k\in \N$, one has
\begin{equation}\label{aniball}
\omega_{a,k} = \frac{2 \pi^{\frac k2}\G(\frac{a+1}2)}{(a+1+k)\G(\frac{a+1+k}2)},\ \ \ \ \ \ \ \sigma_{a,k} = \frac{2 \pi^{\frac k2}\G(\frac{a+1}2)}{\G(\frac{a+1+k}2)}.
\end{equation}
\end{lemma}
\begin{proof}
This can be seen as follows. Cavalieri's principle gives
\[
\omega_{a,k} = \int_{-1}^1 |t|^a \int_{|y| < \sqrt{1-t^2}} dy dt = 2 \omega_k \int_0^1 t^a (1-t^2)^{\frac k2} dt,
\]
where, as customary, we have denoted $\omega_k = \frac{\pi^{\frac k2}}{\G(\frac k2 + 1)}$, the volume of the Euclidean unit ball in $\R^k$.
Keeping in mind that Euler beta function $B(x,y)$, $x, y>0$, can be alternatively expressed as
\[
B(x,y) = 2 \int_0^{1} t^{2x-1}
\big(1-t^2\big)^{y-1} dt,
\]
we immediately recognize that
\[
\omega_{a,k} = \omega_k B(\frac{a+1}2,\frac k2 +1) = \frac{\pi^{\frac k2}}{\G(\frac k2 +1)} \frac{\G(\frac{a+1}2) \G(\frac k2 +1)}{\G(\frac{a+1+k}2 +1)},
\]
which proves the first identity in \eqref{aniball}. To establish the second identity, by \eqref{Za} and the divergence theorem, we have
\begin{equation}\label{oa}
\omega_{a,k} = \frac{1}{a+1+k} \int_{\p B_1}  \sa Z,\nu\da |r|^a d\sigma,
\end{equation}
with $\nu$ denoting the outer unit normal on $\p B_1$. Since by \eqref{eiko} we have $\nu = \nabla \rho$, and by the $1$-homogeneity of $\rho$ we have $\sa Z,\nu\da = Z\rho = \rho = 1$ on $\p B_1$, we infer from \eqref{anivol} and \eqref{oa}
\begin{equation}\label{soa}
\sigma_{a,k} = (a+1+k) \omega_{a,k}.
\end{equation}
This establishes the second identity in \eqref{aniball}.

\end{proof}

By a change of variable, and the translation invariance in the variable ${y}\in \R^k$ of the measure $|r|^a dx$, it is easy to recognise that
\begin{equation}\label{volball}
\int_{B_t(0,{y}_0)} |r|^a dx = \int_{B_t} |r|^a dx= \omega_{a,k}\ t^{a+1+k}.
\end{equation}
On the other hand, Federer's coarea formula and \eqref{volball} gives
\[
\int_0^t \int_{\p B_s(0,{y}_0)} |r|^a d\sigma ds = \omega_{a,k}\ t^{a+1+k}.
\]
Differentiating with respect to $t$ in this identity, and using \eqref{soa}, we conclude that
\begin{equation}\label{sola}
\int_{\p B_t(0,{y}_0)} |r|^a d\sigma = \sigma_{a,k} t^{a+k},\ \ \ \ \ \ \forall {y}_0\in \R^k.
\end{equation}
This leads us to introduce the anisotropic spherical averaging operator
\begin{equation}\label{sa}
M_{a,k}(f,(0,{y}_0),t) = \frac{1}{\sigma_{a,k} t^{a+k}} \int_{\p B_t(0,{y}_0)} f(x) |r|^a d\sigma,
\end{equation}
on functions $f\in C(\R^{k+1})$.
It is clear from \eqref{sola} that
\[
M_{a,k}(f,(0,{y}_0),t)\ \underset{t\to 0^+}{\longrightarrow}\ f(0,{y}_0).
\]
When ${y}_0 = 0$, we simply write $M_{a,k}(f,t)$ instead of $M_{a,k}(f,(0,{y}_0),t)$. The next result generalizes classical properties of the spherical averaging operator.

\begin{lemma}\label{L:sa}
Let $f\in C^1(\R\times \R^k)$. Then for every ${y}_0\in \R^k$ we have
\begin{equation}\label{derM}
\frac{d}{dt} M_{a,k}(f,(0,{y}_0),t) = \frac{1}{\sigma_{a,k} t^{a+1+k}} \int_{\p B_t(0,{y}_0)}  Zf\ |r|^a d\sigma.
\end{equation}
\end{lemma}

\begin{proof}
We begin by observing that for every $t>0$ we have the following alternative representation of the spherical anisotropic averages
\begin{equation}\label{altaa}
M_{a,k}(f,(0,{y}_0),t) = \frac{1}{\sigma_{a,k} t^{a+1+k}} \left[(a+1+k) \int_{B_t(0,{y}_0)} f |r|^a d{x}+ \int_{B_t(0,{y}_0)} Zf |r|^a dx\right].
\end{equation}
If we assume this, by a differentiation in $t$ and the coarea formula, we easily obtain \eqref{derM}.
To see \eqref{altaa}, suppose without restriction that ${y}_0 = 0$. Using the divergence theorem, and the fact that on $\p B_t$ we have $\sa Z,\nu\da = Z\rho = t$, we find
\begin{align*}
M_{a,k}(f,t) & = \frac{1}{\sigma_a t^{a+1+k}} \int_{\p B_t} \sa f |r|^a Z,\nu\da \, d\sigma = \frac{1}{\sigma_a t^{a+1+k}} \int_{B_t} \operatorname{div}(f |r|^a Z) \,dx
\\
& = \frac{1}{\sigma_a t^{a+1+k}} \int_{B_t} f \operatorname{div}(|r|^a Z) \, d{x}+ \frac{1}{\sigma_a t^{a+1+k}} \int_{B_t} Zf |r|^a dx
\\
& = \frac{1}{\sigma_a t^{a+1+k}} \left[(a+1+k) \int_{B_t} f |r|^a d{x}+ \int_{B_t} Zf |r|^a dx\right],
\end{align*}
where in the last equality we have used \eqref{Za}. This gives \eqref{altaa}, thus completing the proof.

\end{proof}

\begin{remark}\label{r:half-space formulas}
The above arguments  were given on domains in $\R\times\R^k$. However, the formulas remain valid for any piece-wise smooth domain  $\Omega\subset (0,+\infty)\times \R^{k}$ with part of its boundary on the singular plane $r=0$ and any function $u\in \mathcal{C}^2(\Om)\cap \mathcal{C}^1(\bar\Om)$. More precisely,  for half-balls $$B^+_t(0,{y}_0)=\{(r,{y})\in B_t(0,{y}_0)\mid r>0 \}\subset \Omega,$$
the used spherical averages can be taken over $\Sigma_t(0,{y}_0)=\p B_t(0,{y}_0)\cap\Omega$  and then use the divergence formula in the domain $ \left( B_t(0,{y}_0)\cap\Omega\right)\subset [0,+\infty)\times \R^{k}$. This is due to the vanishing of the integrals on $\{ r=0\}$ because of the smoothness of $u$ and the weighted surface measure, taking into account that $a>0$.
\end{remark}

We can now prove the strong maximum principle.
\begin{thrm}\label{T:smp}
Let  $\Omega^\star\subset \R^{k+1}$ be such that $J = \Om^\star\cap \{(0,{y})\mid {y}\in \R^k\}\not=\varnothing$. Suppose that $u\in C^2(\Om^\star\setminus J)\cap C^1(\Om^\star)$ be such that in a neighbourhood of any point in $J$ we have $u(r,{y}) = u(-r,{y})$. If $L_a u\geq 0$ in $\Om^\star\setminus J$, then $u$ satisfies the strong maximum principle in $\Om^\star$.
\end{thrm}

\begin{proof}
Let $A = \underset{\Om^\star}{\sup}\ u$. If $A=+ \infty$, there is nothing to prove. Assume therefore that $A<\infty$, and suppose there exist $x_0 = (r_0,{y}_0)\in \Om^\star$ such that $u(x_0) = A$:  we will prove that $u$ is constant in $\Om^\star$. Two cases are possible: either $r_0\not= 0$, or $x_0\in J$. In the former case, suppose to fix ideas that $r_0>0$. Then by the Hopf strong maximum principle for uniformly elliptic equations we would have $u\equiv u(x_0)$ in $\Om$, and therefore in $\Om\cup J$ by the continuity of $u$. Since by assumption $u$ is even in the neighborhood of any point of $J$, we infer that there exists $x_1 = (-r_1,y_1)\in \Om^\star$, with $r_1>0$, where $u$ attains its supremum in $\Om^\star$. Again by Hopf, we infer that $u\equiv u(x_0)$ in all of $\Om^\star$. We are thus left with analyzing the case $r_0 = 0$, i.e. $x_0\in J$. Consider the function $v(x) = A-u(x)$. Clearly $v(x_0) = 0$, $v\ge 0$ and $L_a v \le 0$ in $\Om^\star\setminus J$. By looking at $w(x) = v(r,{y}+{y}_0)$, we can assume without restriction that $x_0 = 0$. By our hypothesis, there exists $\delta>0$, such that $B_\delta\subset \Om^\star$, and for which the function $v(r,{y})$ is even in $r$ in $B_\delta$.
Therefore, for every  $0<t<\delta$ we have from \eqref{derM}
\begin{align*}
\frac{d}{dt} M_{a,k}(v,t) & = \frac{1}{\sigma_{a,k} t^{a+1+k}} \int_{\p B_t}  Zv\ |r|^a  d\sigma = \frac{1}{\sigma_{a,k} t^{a+k}} \int_{B_t}  \operatorname{div}(|r|^a \nabla v)  dx
\\
&  = \frac{1}{\sigma_{a,k} t^{a+k}} \int_{B_t}  L_a v |r|^a  dx \le 0,
\end{align*}
see also \cite[Prop. 12.9]{Gft}.
This implies that for every $\ve<t<\delta$
\[
M_{a,k}(v,\ve) \ge M_{a,k}(v,t).
\]
Letting $\ve\to 0^+$ we find for every $0<t<\delta$
\[
0 = v(0) \ge  M_{a,k}(v,t)\ \Longrightarrow\ 0 \ge \int_{B_\delta} v(x) |r|^a \,d{x}\ge 0.
\]
This shows that $v\equiv 0$ in $B_\delta$, and therefore there exists a point $x_1 = (r_1,{y}_1)\in \Om^\star$, with $r_1\not= 0$, where $u(x_1) = A$. Again by the Hopf strong maximum principle, we infer that $u\equiv A$ in $\Om^\star$.

\end{proof}

\begin{remark}
We recall that $a\ge 0$. When $a>0$, in the above application of the divergence theorem one has to be a bit careful and split the integral on $B_t$ into two integrals on the region $B^{\ve}_t = B_t \cap \{|r|>\ve\}$, and then let $\ve\to 0^+$. The boundary integrals on $B_t \cap \{|r| = \ve\}$ converge to zero as $\ve\to 0^+$ because $a>0$.
 \end{remark}

\section{The Serrin type problem}\label{S:main}

The main goal of this section is to prove Theorem \ref{t:serrin}.
We begin with a result which allows an (even) $C^1$ continuation of a solution to the equation $L_a = - 1$ across the singular plane $r=0$.

\begin{prop}\label{p: normal der vanish}
Let $a\not=0$ and $D=(0,R)\times  (b_0, b_1)^k$. Suppose that $u\in C^2(D)\cap C^1(\bar D)$ satisfy $L_a u= g$, where $g\in L^1 (D)$. Then, for any ${y}\in  (b_0,b_1)^k$ we have
\[
u_r(0,{y})=\lim_{r\rightarrow 0+} u_r(r,{y})=0.
\]
\end{prop}

\begin{proof}
Since the equation is translation invariant in the variable $y$, it is enough to prove the claim in the case  $D=(0,R)\times (0,b)^k$.
Let $0<\varepsilon<b$, $0<t<R$, and denote
\[
D'=(0,t)\times Q' _k= (0,t)\times (\varepsilon,b)^k.
\]
We also indicate with $Q^i _{k-1}$ the $(k-1)$-dimensional parallelepiped $(\varepsilon,b)^{k-1}$, in which the variable $y_i$ is missing, and we let
${d\mathcal{L}_{k-1}}$ be the $(k-1)$-dimensional Lebesgue measure.
First, we record the following identity which is essentially contained  in \cite[p. 374]{Walter57} in the homogeneous case  $g=0$, see also \cite[(4.2)]{Fox59},  where the identity is proven for ultrahyperbolic equations, and \cite[Lemma p. 131]{Young72}, where an even more general operator is considered.
Let $F=u_r^2 - \sum_{i=1}^ku^2_{{y}_i}\in C(\bar D)$. A calculation shows that we have
\begin{equation}\label{e:walter's}
2 r u_r g = 2r u_r\ L_a u = \left[ r\, F \right]_r + 2r\sum_{i=1}^k \left(u_{{y}_i}u_r\right) _{{y}_i} -  F +2a u_r^2.
\end{equation}
Integrating \eqref{e:walter's} over  $D'$ and using the divergence theorem we obtain  the identity
\begin{align*}
2 \int_0^t \int_{Q' _{k-1}}  r u_r g {dy} dr  & = t\int_{Q'_k}F(t,y) dy + 2\sum_{i=1}^k \int_0^t \int_{Q^i _{k-1}}r (u_ru_{{y}_i}) \big|_{{y}_i=\varepsilon} ^{{y}_i=b} {d\mathcal{L}_{k-1}}dr
\\
& -\int_0^t \int_{Q' _{k}} \left( F- 2 a u_r^2\right) {dy} dr.
\end{align*}
Dividing by $t$ in this identity,  we find
\begin{equation}\label{e:int normal 1}
2 \int_0^t \frac {r}{t} \int_{Q' _{k-1}}  u_r g  {dy} dr= I_1+I_2+I_3,
\end{equation}
where
\[
I_1=\int_{Q'_k}F(t,y) dy,
\quad I_2=2\sum_{i=1}^k \int_0^t \frac {r}{t}  \int_{Q^i _{k-1}} (u_ru_{{y}_i}) \big|_{{y}_i=\varepsilon} ^{{y}_i=b}  {d\mathcal{L}_{k-1}}dr,
\]
and
\[
I_3=-\frac{1}{t} \int_0^t \int_{Q' _{k}} \left( F- 2 a u_r^2\right){dy}dr.
\]
Since $0<\frac {r}{t} <1$ and $u_r g\in L^1 (D)$  by the assumptions, we have
\[
\left|\int_0^t \frac {r}{t} \int_{Q' _{k-1}}  u_r g  {dy} dr\right| \le \int_0^t \int_{Q' _{k-1}}  |u_r| |g|  {dy} dr\ \underset{t\rightarrow 0+}{\longrightarrow}\ 0.
\]
Similarly, we obtain
\[
\lim_{t\rightarrow 0+} I_2=0.
\]
Furthermore, since $F\in C(\bar D)$, one has
\[
\lim_{t\rightarrow 0+} I_1=\int_{Q'_k}F(0,y) dy.
\]
Finally, the fundamental theorem of calculus gives
\[
\lim_{t\rightarrow 0+} I_3=-\lim_{t\rightarrow 0+}  \frac {1}{t} \int_0^{t}\int_{Q' _{k}} F- 2 a u_r^2 \ {dy}d\tau= - \int_{Q' _{k}} F(0,y) {dy}
- 2 a \int_{Q' _{k}}  u_r^2(0,y)  {dy}.
\]
In conclusion, since $a\not=0$, we have obtained on $Q' _k=  (\varepsilon,b)^k$ the identity
\[
0=\lim_{t\rightarrow 0+} \left( I_1 +I_2+I_3\right)=\int_{Q' _{k}}  u_r^2(0,y) {dy}.
\]
Thus, for every $0<\varepsilon<b$, we have the vanishing of  $ u_r^2(0,y)\equiv 0$ on  $Q' _k=  (\varepsilon,b)^k$. Since $u\in C^1(\bar D)$ we obtain that $ u_r(0,y)\equiv 0$.

\end{proof}

Our next objective is to prove Theorem \ref{t:P integral for serrin}, an integral constraint for the function $P_a$ defined by \eqref{Pa}. In the sequel we assume that
$\Omega^\star\subset \R^{k+1}$ be a smooth, piecewise $C^1$ connected domain, symmetric with respect to the $r=0$, and as before let $x=(r,y)\in \R^{k+1}$. Recall that we are denoting $\Omega$  the part of $\Omega^\star$ in the half-space $\R^{k+1}_+$. Thus, $\Omega^\star$ is the reflected double of $\Omega$ with respect to the hyperplane $r=0$ in $\R^{k+1}$. Let  $\Sigma=\partial\Omega^\star\cap\{r> 0\}$.  In particular, the boundary $\partial \Sigma $ of $\Sigma$ is a smooth $(k-1)$-dimensional surface in the plane $\{r=0\}$. Finally, $\bomega$ is the disjoint union of $\Sigma$, the flat part $\Sigma_0=\bar\Omega^\star\cap \{r>0\}$, and the  $(k-1)$-dimensional surface $\Sigma'=\partial \Sigma= \partial \Sigma_0 $ in the plane $\{r=0\}$,
\begin{equation}\label{e:bdry notation}
\bomega=\Sigma\cup\Sigma_0\cup \Sigma'.
\end{equation}

In the statement of the following proposition $\gr$ denotes the gradient with respect to the variable $x\in \R^{k+1}$ and $\sa \cdot,\cdot \da$ denotes the Euclidean scalar product. Its proof is inspired to that of Weinberger's ingenious use of the identity first established by Rellich in \cite{Relid}.

\begin{prop}\label{p:Pohozaev}
Suppose $u\in C^2(\Om)$  is a solution of  $L_a u = - 1$ in $\Om$. If $Z$ is the radial vector field \eqref{e:radial field}, then the following identity holds true,
\begin{equation}\label{e:pohozaev}
\divg \left ( r^a\frac {|\gr u|^2}{2}Z  - r^a Zu \nabla u- r^au Z\right)  = [(a+1+k)-2)] r^a\frac {|\gr u|^2}{2}  - (a+1+k) r^a u.
\end{equation}
\end{prop}

\begin{proof}
We begin with the key observation that, if for convenience we indicate $\nabla = (D_1,...,D_{k+1})$, then the commutator $[D_j,Z]$ satisfies the identity
\[
[D_j,Z] = D_j,\ \ \ \ \ j=1,...,k+1.
\]
This has the following direct consequence
\begin{equation}\label{one}
Z(\frac{|\nabla u|^2}2) = \sa\nabla(Zu),\nabla u \da - |\nabla u|^2.
\end{equation}
Furthermore, we have
\begin{equation}\label{two}
\operatorname{div}(r^a Zu \nabla u)  = Zu \operatorname{div}(r^a \nabla u) + r^a \sa\nabla(Zu),\nabla u \da
 = - r^a Zu + r^a Z(\frac{|\nabla u|^2}2) + r^a |\nabla u|^2,
\end{equation}
where in the last equality we have used \eqref{one} and the fact that, since $L_a u = r^{-a}\divg \left(r^a\gr u \right)$, the equation $L_a u = - 1$  is equivalent to $\divg \left(r^a\gr u \right)=-r^a$. Since by \eqref{Za} we have  $\divg (r^a Z) = (a+1+k)r^a$, we conclude from \eqref{two}
\begin{align*}
& \divg \left ( r^a\frac {|\gr u|^2}{2}Z  - r^a Zu \nabla u- r^au Z\right) = (a+1+k) r^a \frac{|\gr u|^2}2  + r^a Z(\frac {|\gr u|^2}{2})
\\
& + r^a Zu - r^a Z(\frac{|\nabla u|^2}2) - r^a |\nabla u|^2 - (a+1+k) u - r^a Zu
\\
& = (a-1+k) r^a \frac{|\gr u|^2}2  - (a+1+k) r^a u,
\end{align*}
which completes the proof of \eqref{e:pohozaev}.

\end{proof}

Our next goal is to establish the following integral identity for the function $P_a$ which has been introduced in \eqref{Pa}.

\begin{thrm}\label{t:P integral for serrin}
If $u\in C^2(\Om)\cap C^1(\bar\Om)$ is a solution to
 \begin{equation}\label{e:BV for P integral for serrin}
L_a u = - 1\quad \text{in}\ \Om\qquad\quad
u _{\big|\Sigma}=0,\qquad |\nabla u|_{\big|\Sigma} = c,
\end{equation}
then the following identity holds true
\begin{equation}\label{e:P fn integral identity}
 \IO \left ( P_a(x)-c^2\right) \vola=0.
\end{equation}
\end{thrm}

\begin{comment}
 The proof of Theorem \ref{t:P integral for serrin} will follow from the following claims which involve the already  function $P_a$,  cf. \eqref{Pa},
Note that we have: (i) $L_aP_a\geq 0 $ in $\Omega$ by \eqref{wow},  (ii) the maximum principle, cf. Theorem \ref{T:smp}, and (iii) the trivial identity $\quad P_a\vert_{\Sigma}=c^2$, which follows from the assumed boundary conditions.
 Our next goal is to prove:
\begin{enumerate}
\item The normal derivative satisfies $u_{x_0} = 0$ on ${x_0}=0$;
\item A Pohozaev identity, which will be used to see $\int_{\Omega^*_{+}} \left ( P_a(x)-c^2\right) x_0^a dx=0;$
\item The fact that   $$u=\frac {R^2-x_0^2- |x' - x'_0|^2}{2(a+k+1)}$$ and  $\Omega^*$ is the ball
$B_R((0,z'))=\{ (x_0,x')\mid x_0^2+ |x' - z'|^2<R^2\}$ with $$\Omega^*_{+}=\{ (x_0,x')\mid x_0^2+ |x' - z'|^2<R^2, \ x_0>0\}.$$
\end{enumerate}
\end{comment}

\begin{proof}
Let ${\nu}$ be the  outer unit normal to $\Omega$ and $d\sigma$ be the $k-$dimensional surface measure on $\bomega$.  Since the radial vector field $Z$ is smooth and $a\ge 0$ we have  $\divg \left( r^a Z\right)=(a+1+k)r^a\in C\left([0,+\infty)\times\R^k \right)$. Hence, the divergence theorem gives
\begin{equation}\label{e:div prop 3}
(a+1+k)\int_\Omega \vola = \int_\bomega \sa Z,{\nu}\da \surfa.
\end{equation}
On the other hand, the divergence theorem and the equation $u\divg \left(r^a\gr u \right)=-r^au$, give
\begin{equation}\label{e:div prop 1}
-\int_\Omega  |\gr u|^2 \vola  + \int_\bomega u u_\nu \surfa = -\int_\Omega u \vola.
\end{equation}
We observe now that
\[
\int_\bomega u u_\nu \surfa = \int_{\Sigma} u u_\nu \surfa + \int_{\Sigma_0} u u_\nu \surfa = 0.
\]
The former integral in the right-hand side vanishes by the boundary condition $u = 0$ on $\Sigma$ in \eqref{e:BV for P integral for serrin}, the latter does as well because $r^a = 0$ on $\Sigma_0$.
We thus obtain from \eqref{e:div prop 1}
\begin{equation}\label{e:div prop 2}
\int_\Omega  |\gr u|^2 \vola   = \int_\Omega u \vola.
\end{equation}

The identity \eqref{e:pohozaev} and the divergence theorem give
\begin{multline}\label{e:pohozaev 1}
\int_\bomega \left ( \frac {|\gr u|^2}{2}\sa Z,{\nu}\da - u_\nu Zu    - u \sa Z,{\nu}\da\right) \surfa \\
 =\ (a-1+k)\int_\Om  \frac {|\gr u|^2}{2}\vola  - (a+1+k)\int_\Om u \vola.
\end{multline}
The condition $u = 0$ on $\Sigma$ gives $\nabla u = u_\nu \nu$, and therefore we have on $\Sigma$
\[
u_\nu Zu = u_\nu\sa Z,\gr u\da = u_\nu^2 \sa Z,\nu\da = c^2 \sa Z,\nu \da,
\]
by the constant Neumann boundary condition.
Using this identity in \eqref{e:pohozaev 1}, together with  the Dirichlet boundary condition and the vanishing of $r^a$ on the flat part of the boundary $\Sigma_0$, we obtain
\begin{equation}\label{e:pohozaev 1a}
-\frac {c^2}{2}\int_\bomega  \sa Z,{\nu}\da \surfa
 = \frac{a-1+k}{2}\int_\Om  {|\gr u|^2}\vola  - (a+1+k)\int_\Om u \vola.
\end{equation}
Combining \eqref{e:pohozaev 1a} with \eqref{e:div prop 3}, we find
\[
\frac{2-(a+1+k)}{2}\int_\Om  {|\gr u|^2}\vola  + (a+1+k)\int_\Om u \vola = \frac {c^2}{2}(a+1+k)\int_\Omega \vola.
\]
If we now use \eqref{e:div prop 2} in the latter equation, we obtain
\begin{equation}\label{e:pohozaev 4}
((a+1+k)+2)\int_\Om u\vola = c^2(a+1+k)\int_\Om\vola,
\end{equation}
which, after dividing by $a+1+k$, again by \eqref{e:div prop 2} can be rearranged as
\[
\int_\Om \left[|\nabla u|^2 + \frac{2}{a+1+k} u\right]\vola = c^2\int_\Om\vola.
\]
Recalling \eqref{Pa}, we have reached the desired conclusion \eqref{e:P fn integral identity}.

\end{proof}

We are finally ready to provide the proof of our main result.

\begin{proof}[Proof of Theorem \ref{t:serrin}]
We consider only the case $a>0$ since otherwise we are in the framework of Serrin's result.
We begin by noting that $u$ must be an even function in the $r$ variable. Indeed, taking into account Proposition \ref{p: normal der vanish} and the facts that both the equation  and the domain are  invariant under reflections with respect to the plane $\{r=0\}$, we have that $$v(r,y)=u(r,y)-u(-r, y)\in C^2(\Om^\star\setminus \{r=0\})\cap C^1(\bar\Om^\star)$$ and
\[
L_a \, v=0 \ \text{in} \ \Om^\star\setminus \{r=0\}),\qquad  v = 0\quad  \text{on}\ \p \Om^\star.
\]
By Theorem \ref{T:smp} we conclude that $v\equiv 0$, i.e., $u$ is even with respect to $r$. Once we know this, we can appeal to \cite[Theorem 1.1]{STV}, and infer that, in fact, $u\in C^\infty(\Om^\star)$, and therefore $P_a\in C^\infty(\Om^\star)$. Furthermore, note that the evenness in $r$ of $u$ implies that $P_a(-r,y) = P_a(r,y)$ in $\Om^\star$. We infer that $P_a$ satisfies the hypothesis of the strong maximum principle in Theorem \ref{T:smp}, since
$L_a P_a \ge 0$ in $\Om^\star\setminus J$ by \eqref{wow}. Noting that the assumed boundary conditions for $u$ show that  $P_a=c^2$ on $\p\Om^\star$, we infer that either $P_a< c^2$ or $P_a\equiv c^2$   in ${\Om^\star}$. The former possibility, however, is ruled out by applying  Theorem \ref{t:P integral for serrin} to each of the sides with respect to the plane $r=0$ of $\Om^\star$. By the evenness of $P_a$ with respect to $r$, we obtain the integral identity
\begin{equation}\label{e:P fn integral identity Om*}
 \int_{\Om^\star} \left (P_a(x)-c^2\right) |r|^a dx =0.
\end{equation}
We conclude that it must be $P_a\equiv c^2$ in $\Om^\star$. This implies, in particular, that $L_a P_a \equiv 0$ in $\Om^\star$.
The end of the proof now follows from the case of equality in \eqref{wow} in Theorem \ref{T:P}.

\end{proof}

%%%%%%%%%%%%%%%%%%%%%%%%%%%%%%%%%%%%%%%%%%%%%%%%%%%%%%%%%%

\end{document}